\documentclass[12pt]{amsart}
\usepackage{graphicx}
\usepackage{amssymb}
\usepackage{setspace}

\usepackage[all]{xy}

\newcommand{\C}[1]{{\mathcal#1}} 

\theoremstyle{plain}

\theoremstyle{definition}

\theoremstyle{definition}
\newtheorem{example}{Example}[section]

\theoremstyle{remark}
\newtheorem{rem}{Remark}[section]

\usepackage{fancyhdr}
\begin{document}
\title{On Periods: from Global to Local} 
\author{Lucian M. Ionescu}
\address{Department of Mathematics, Illinois State University, IL 61790-4520}
\email{lmiones@ilstu.edu}
\date{June, 2018} 

\begin{abstract}
Complex periods are algebraic integrals over complex algebraic domains,
also appearing as Feynman integrals and multiple zeta values.
The Grothendieck-de Rham period isomorphisms
for p-adic algebraic varieties defined via Monski-Washnitzer cohomology, is briefly reviewed.

The relation to various p-adic analogues of periods are considered,
and their relation to Buium-Manin arithmetic differential equations.
\end{abstract}

\maketitle
\setcounter{tocdepth}{3} 
\tableofcontents

\section{Introduction}
In this article we discuss periods and their applications, 
as a continuation of \cite{Ionescu-Sumitro}, focusing on the relation between global periods in characteristic zero,
and their local counterparts.

The main goal of this research is to question the ``stability'' of the connection between scattering amplitudes and periods \cite{Schnetz:QuantumPeriods,Brown:FeynmanIntegrals,Brown:ICMP}, when passing from global to local, by using the analogy between Veneziano amplitudes and Jacobi sums,
is addressed in a follow up article \cite{LI:p-adicFrobenius},
which also adopts a Deformation Theory point of view when introducing p-adic numbers.
It is expected to provide some feedback on the Feynman amplitudes and Multiple Zeta Values correspondence \cite{QuantaMagazine,LI:Periods-FI-JS-Talk}.

\vspace{.1in}
Periods are values of algebraic integrals, extending the field of algebraic numbers.
Non-trivial examples are Feynman amplitudes from 
experimentally ``dirty-gritty'' Quantum Field Theory \cite{Schnetz:QuantumPeriods},
yet which happen to be also linear combinations of multiple zeta values from ``pure''
Number Theory \cite{QuantaMagazine}. 
That Mathematics is unreasonably effective, we know; but to the point of
starting to reconsider Plato's thesis that reality is a mirror of the world of 
(mathematical) ideas?! So, {\em Number}, (once categorified) does rule the (Quantum) Universe after all ...

After reviewing the idea and concept of period, the article explores the connection with {\em quantization functors}, i.e. representations of (generalized) categories of cobordisms, 
as a perhaps more physical route than that of abstract motives.

At a more concrete level,
the power of analogy \cite{Weil-analogy} between Veneziano amplitude as a 
String Theory analogue of Feynman amplitude,
and Jacobi sum in finite characteristic \cite{Ireland-Rosen} 
(an finite characteristic analog of Euler beta function),
is used to investigate a possible global-to-local correspondence for periods (factorization or reduction of cohomology):
$$\xymatrix@R=.2pc{
Veneziano\ Amplitude: & & Jacobi\ Sum: \\ A(a,a')=\frac{\Gamma(\alpha)\Gamma(\beta)}{\Gamma(\alpha+\beta)} 
& \quad \leftrightarrow \quad & J(c,c')=\frac{g(c)g(c')}{g(cc')}
}
$$
where $\alpha=-1+(k_1+k_2)^2, \beta=-1+(k_3+k_4)^2$ relate the in/out momenta of the interacting strings, 
and $c,c':F_p^\times->C^\times$ are multiplicative characters of the finite field $F_p$.

The first measures the correlation (interaction amplitude) of two strings, with momenta expressed in Mandelstam's variables,
while the second measures, the ``intersection correlation'' between two multiplicative subgroups (e.g. squares and cubes),
yielding the correction term (``defect'' $a_p$) for the number of points $N_p$
(like a constructive or destructive amplitude for the ``volume integral'') 
of a finite Riemann surface $C(F_p)$, over a finite field $F_p$.
This connection has been studied, as for example in \cite{Kholodenko}.

The Local-to-Global Principle could be used even informally, via an analogy with the algebraic number theory case, to guide more experienced investigators deal with the global case of periods, 
Feynman Diagrams and Mirror Symmetry.
In the other direction, it can be used to guide the development of p-adic String Theory, beyond a mere formal replacement of real (complex) numbers by p-adic numbers.

The interplay between Galois symmetries and periods (Feynman amplitudes) \cite{Brown:ICMP}, 
will be investigated in the framework of Noether's Theorems, 
connecting conserved quantities (and e.g. unitarity as conservation of probability), and symmetries of the system.

The article is organized as follows.
We review the basic ideas regarding periods in \S \ref{S:Periods-AG}, 
starting from their simple introduction as algebraic integrals, followed by a cohomological interpretation.
Remarks on periods, motives and Galois group are followed by considering p-adic periods, in connection with Buium calculus.

Further considerations are postponed for a Deformation Theory approach to p-adic numbers \cite{LI:p-adicFrobenius}, and an investigation of a connection between Grothendieck's algebraic de Rham cohomology, and the discrete analog of de Rham cohomology of the present author \cite{LI:DiscreteDeRham}, as well as possible connections with the discrete periods of \cite{Mercat:DiscretePeriods}.

\section{Periods: from integrals to cohomology classes}
The arithmetic notion of period refers essentially 
to the value of rational integrals over rational domains \cite{KZ}.
For example the ubiquitous ``Euclidean circle-radius ratio'' $\pi=\int_{[-1,1]}dx/\sqrt{1-x^2}$, 
residues like $2\pi i = \int dz/z$ or path integrals as $\log(n)=\int_1^n dx/x$ \cite{K},
are {\em numeric periods}.

The representation of a period as an integral is not unique. 
When placed in the context of de Rham isomorphism of a compact manifold,
or algebraic variety as in the early work by Grothendieck \cite{Grothendieck-deRham},
they are matrix coefficients of the corresponding integration pairing \cite{Muller-Stach}.
The resulting isomorphism between de Rham cohomology and singular cohomology:
$$de \ Rham\ Theorem: \quad 
H_{dR}^*(X,D)\otimes_Q C \quad \overset{\cong}{\to}\quad H_{sing}^*(X,D)\otimes_Q C,$$
is called the {\em period isomorphism}.

\subsection{Periods of Algebraic Varieties} \label{S:Periods-AG}
Specifically, from the algebraic geometry point of view, 
an {\em numeric period} $p$ is represented by a 
quadruple consisting of a smooth algebraic variety $X$ over $Q$,
of dimension $d$,
$\omega$ is a regular algebraic d-form, $D$ a normal crossing divisor 
and $\gamma$ is singular chain on $X(C)$, with boundary on $D(C)$:
$$(X,D,\omega,\gamma) \quad \mapsto \quad p=\int_\gamma \omega.$$
Fixing $X$ and such a normal divisor $D$, 
choosing a rational basis in both cohomology groups allows to 
represent the above {\em period isomorphism} as a {\em period matrix}.

Of course, there are elementary transformations on such quadruples (linearity,
change of variables and Stokes formula),
which leave the corresponding period unchanged \cite{KZ}, p.31.
Whether the {\em effective periods}, i.e. 
equivalence classes of quadruples modulo the elementary moves, 
correspond isomorphically to the numeric periods, 
is the content of the corresponding Kontsevich Conjecture.

\vspace{.2in}
Since the history of periods and period domains goes back to the very beginning of algebraic geometry \cite{CG}, we will proceed with two such elementary examples:
the Riemann Sphere (genus zero), the case of elliptic curves (genus une).

\begin{example}
With $X=P^1-\{0,\infty\}$, $D=\empty$, $\omega=dz/z$ and $\gamma=S^1$ the unit circle,
we find $2\pi i$ as the (only) netry of the period matrix of the 
period isomorphism $H^1(X;C)\to H_1(X;C)$.

If we change the divisor to $D=\{1,n\}$ and take $\gamma=[1,n]$, then
the numeric period $\log(n)$ becomes one of the periods of
$H^1(X,D)$.
\end{example}

\begin{example}\label{Example:EC}
Given an elliptic curve $X:y^2=x^3+ax+b$ 
with canonical homology basis $\gamma_1,\gamma_2$ \cite{RS} 
and differential form $\omega=dx/y$, the period matrix (vector) is (\cite{CG}, p.1418):
$$(A,B)=\left( \int_{\gamma_1}\omega, \int_{\gamma_2} \omega\right).$$
It is customary to go to the fraction field (of periods)
and divide by $A$ to get $(1,\tau)$ with the {\em normalized $B$ period}
denoted $\tau=1+it$, 
having positive imaginary part\footnote{Due to the fact that 
$i\int_X \omega\cup \bar{\omega}>0$ \cite{CG}.}
and constituting an invariant of the elliptic curve \cite{Carlson-Stach-PD}, p.9.

For example with $\lambda=-1$, the elliptic curve $E:y^2=(x-1)x(x+1)$
has invariant $\tau=i$ and $E=C/Z\oplus Zi$ has an additional isomorphism 
of order $4$ (complex multiplication) \cite{Chowla-Selberg}.

Such normalized periods provide a simple example of {\em period domain},
here the upper-half plane $\C{H}$.

A change of the homology basis by a unimodular transformation in 
$\Gamma=SL_2(Z)$ corresponds to a fractional transformation relating 
the corresponding two points of the period domain.
Thus the moduli space of genus one Riemann surfaces corresponds to
the quotient $\Gamma / \C{H}$.

Additional details and examples can be found in \cite{Carlson-Stach-PD}, Ch.1.
\end{example}

Before we move on, a crude physics interpretation of periods will be later useful.
\begin{rem}
Think of an elliptic curve with a 1-form $\omega$ as a 2D-universe with a
flow, or perhaps a world sheet of a string, with a given capacity for 
propagating action (Quantum computing: duplex channel).
In a conformal metric, for convenience of relating to metric picture,
the corresponding vector field will represent a free propagation,
with certain circulation and flux 
(harmonic dual pair: streamlines and equipotential lines).
The two periods then measure these two: circulation and flux.
\end{rem}

\subsection{Families of periods}
Continuing the discussion of the above case of Riemann surfaces of genus one, 
one often has a family of such elliptic curves $E_t$ depending holomorphically 
on a parameter $t$.
The resulting map $t\mapsto \tau(t)$ is the {\em period map}.

For example, if the base space is the Riemann sphere,
one finds a globally defined map $S\mapsto \Gamma/\C{H}$ \cite{CG}, p.1418.

%
%
\subsection{Interpretation}
Comparing periods and algebraic numbers is probably the first thing to do, 
before developing a theory of periods.

\subsubsection{On algebraic numbers}\label{S:AlgebraicNumbers}
Algebraic numbers (over $Q$) extend rational numbers via extensions $Q[x]/(f(x))$.
If focusing on integers, and contenting to systematically view field extensions 
as fields of fractions (as long as we stay within the commutative world),
then we may choose to interpret these algebraic extensions geometrically,
as lattices \cite{IonescuMina}, and algebraically as {\em group representations}.
For example, $Q(i)$ is the fraction field of its ring of algebraic integers $Z[i]$,
which in turn is the group ring of its group of units $U=<i>$,
which is a subgroup of the rotation group of the {\em rational plane} $ZxZ$
\footnote{... in the spirit of the geometric interpretation of complex numbers, starting with Argand, Gauss, followed by Riemann and perhaps the modern CFT and String Theory developments.}

At this stage, finite fields (finite characteristic) can be constructed via quotients 
of lattices of algebraic integers.
For example $F_5\cong Z[i]/(2+i), F_{3^2}=Z[i]/(3)$ etc. 
\footnote{See \cite{IonescuMina} for developments of this direction of reasoning.}.

Then, increasing the dimension from the above $D=0$ case, i.e. the number of variables, 
we obtain algebraic varieties $X:Z[x_1,...,x_n]/<f_1,...,f_k>$, 
suited for a geometric study via homology and cohomology.
As a reach class of examples, we have the Riemann surfaces
$RS:y^2=f(x)$, with de Rham / Dolbeault cohomology over the complex numbers.

\subsubsection{Cohomology pairings and their coefficients}
Integrals, as a non-degenerate pairing, e.g. de Rham isomorphism, 
on the other hand go beyond the arithmetic realm.
Nevertheless periods should probably be beter compared with algebraic integers,
since thw inverse of a period is not necessarily a period: 
$1/\pi$ is conjecturally not a period.

Another important aspect is that other important numbers like $e$ and Euler's 
constant $\gamma$, conjecturally are not periods; why? what lies still beyond 
periods, but before the ``transcendental junk'' like Louisville's number and such!?
Is there a numerical shadow of a Lie correspondence, and hence
exponential periods like $e^{period}$ form another important class?

\subsubsection{Speculative remarks}
Speaking of shadows, the cohomology theory lies ``above'' the periods themselves,
reminiscent of {\em categorification}, 
a process relating numbers and algebraic structures
(e.g. Grothendieck ring etc.).

But how to isolate the ``cohomology theory'' from the actual implementation,
based on a specific manifold? 
Is there such a ``thing'' as isomorphism class of the functor representing the 
respective cohomology theory? And how are its various matrices,
corresponding to its values, related?

If algebraic numbers can be viewed in fact as representations 
corresponding to their multiplicative structure, as suggested above,
then periods should probably relate to representations of groupoids, maybe?

This leads to the ``3-rd level'' of abstraction, beyond arithmetic and algebraic-geometric.

%
%
\subsection{Periods, Motives and Galois Group}
\cite{Wiki-motives} ``The theory of motives was originally conjectured as an attempt to unify a rapidly multiplying array of cohomology theories, including Betti cohomology, de Rham cohomology, l-adic cohomology, and crystalline cohomology.''

In a more abstract direction, following the work of Grothendieck on pure and mixed motives,
in 1990s the work of M. Nori starting from directed graphs encoding ``equivalence moves'' between periods (with a certain similarity to Redemeister's moves and theorem on homotopy classes of knots), let to a degree of abstraction which seams not to serve our purpose here,
to understand the ``knots, braids and links'' of {\em real} ``elementary'' particles
in decays and collisions.

Now on one hand, Feynman integrals, also periods, are ``closer'' to Chen's theory of iterated integrals, which forms a homotopy theory analog of de Rham cohomology \cite{Chen,Hains}.
This explains the ``coincidence'' with (linear combinations of) periods arising from Number Theory, e.g. multiple zeta values, which can also be expressed as iterated integrals.
This ``non-commutative side'' of {\em homotopical motives}, is probably better suited to
be approached from the (``Cosmic'') Galois action viewpoint \cite{Brown:Galois}.

On one the other hand,
the de Rham {\em cohomology} pairing framework for understanding and generalizing periods
is reasonably close/similar to the {\em representation theory} viewpoint for algebraic integers.

Then what is the connection between the two directions? 
It seams that generalizing the idea of Galois action on roots of polynomials (representations
point of view), allows to view periods as having infinite orbits under bigger analogues of
``Galois groups'' \cite{Brown:Galois}. 
For example $\pi$ can be viewed as associated to a 1-dimensional representation of a group
\cite{Brown:FeynmanIntegrals}, p.11.
This direction provides a framework for studying amplitudes,
perhaps not unrelated to the cornerstone idea in high energy physics
that ``elementary particles'' are associated to irreducible representations,
but definitely to be pursued by ``Mathematicians only!''.

Alternatively, {\em motives} as universal sumands of cohomology theories, 
allows to connect with the direct approach to periods via the period isomorphism.

This perhaps allows to connect the global and the local.
Indeed, since {\em Weil cohomologies} are such universal sumands, one should be able 
identify what are the natural analog of the concept of period in finite characteristic
\cite{Weil-Conjectures}.

We will confine to the simplest nontrivial case of elliptic curves (Example \ref{Example:EC}),
and investigate in what follows the reduction modulo a prime,
in the context of {\em Ramification Theory}, 
within the {\em Algebraic Number Theory framework},
suited for the algebraic integers and representation viewpoint,
and which is ``close enough'' to 
Algebraic Geometry framework mentioned in \S\ref{S:Periods-AG}.

\section{p-Adic Periods}
The above periods in characteristic zero correspond to the ``real world'' 
of the ``prime at infinity'' (arguably on both accounts: \cite{Real-fish}
\footnote{Real numbers result in completing the rationals the ``other way''
then the direction of the carryover 2-cocycle!}).

But what about the periods of p-adic analysis?
How are these defined, and how do they relate to p-adic analogues like p-adic gamma function, Gauss sums and Jacobi sums?
Of course, beyond the natural motivation to generalize and applications to Number Theory, such a study would share light on p-adic String Theory and CFT, via their connection to Veneziano amplitudes and other such 
iterated integrals on moduli spaces of punctured Riemann spheres \cite{Brown:ModuliSpaces}.

\subsection{p-Adic De Rham Cohomology}
(Algebraic/Geometric) Number Theory in finite characteristic 
may be thought of as the ``infinitesimal/linear analysis''
of p-adic analysis\footnote{It is rather {\em Deformation Theory}, 
as it will be argued elsewhere \cite{LI:p-adicFrobenius}.},
and algebraic de Rham cohomology of a variety does not reduce ``nicely'',
requiring a lift to characteristic zero of p-adic number fields,
called Monski-Washnitzer cohomology \cite{Hartog}, p.27 (See also \cite{Kedlaya} and references therein).

Briefly, if $X$ is an algebraic variety over $Q$ and $A=Q[x,y]/f(x,y)$ its coordinate ring,
one considers the ``overconvergent'' subalgebra $A^\dagger$ of its lift to p-adic numbers (loc. cit.),
in order to have exact forms closed under p-adic completion.
Then the {\em Monski-Washnitzer cohomology of $X$},
also called here the {\em p-adic algebraic de Rham cohomology of $X$} is:
$$H^i_{MW}(A):= H^i_{dR}(A^\dagger)\otimes_{Z_q} Q_q,$$
where $q=p^n$, $Z_q$ is the p-adic integers of the $n-th$ degree unramified extension of $Z_p$,
and $Q_q$ its field of fractions.

Constructing a lift to $Z_p$ of Frobenius $x\mapsto x^q$, as a ring endomorphism turns out to be difficult,
being equivalent to specifying a p-derivation in the sense of Buium-Manin \cite{Buium-Manin}:
$$\phi_p(x)=Frob_p(x)+p\delta_p(x), \quad \delta(x+y)=\delta(x)+\delta(y)+C_p(x,y),$$
$$C_p(x,y)=[x^p+y^p-(x+y)^p]/p\ \in \ Z[x,y].$$
Note at this stage that, if concerned with its action on MW-cohomology, then
one may relax the endomorphism requirement, and examples
like $\phi_1(x)=x^p$ or $\phi_2(x)=x^p+px$ will do \cite{Hartog}, Ch.3, p.24.

\subsection{p-adic Period Isomorphism}
The p-adic analog of period isomorphism is defined via Hodge Theory and Hodge isomorphism.

Since the presentation may benefit from a deformation Theory approach to p-adic numbers, as sketched in \cite{LI:p-adicFrobenius}, it will continued in loc. cit.
Indeed, the understanding of the p-adic periods would benefits from a comparison of the MW-cohomology with Hochschild cohomology, 
and corresponding period isomorphism. 

\section{Conclusions}
The subject of {\em Periods}, although with a long history, has become essential for understanding the deep mysteries underlying the ``coincidence'' between quantum physics scattering amplitudes and number theoretical special values, like MZVs.

Since the Grothendieck's algebraic de Rham cohomology is instrumental in studying the global (conceptual) aspects of periods, we will also point to its connection to classical de Rham cohomology, via the discrete version of de Rham cohomology, as defined for finite abelian groups \cite{LI:DiscreteDeRham}, which will be addressed elsewhere, in connection with p-adic periods.

Finally, the analogy between Euler's integrals, beta and gamma functions, with the Jacobi and Gauss sums, is used to question the ``stability'' of the (Veneziano) amplitude-periods connection, when passing from global to local.

Understanding the connections with Venatiano amplitudes and Jacobi sums require an understanding of the ``discrete case'', of finite characteristic.
A parallel with characteristic zero can be achieved via a Deformation Theory viewpoint, and can be found in \cite{LI:p-adicFrobenius}.

\section*{Acknowledgements}
The author would like to thank IHES for the excellent conditions making possible the continuation of author's research on a theme of Professor Kontsevich, under such a highly valuable guidance. 


\end{document}